\documentclass[review]{elsarticle}
\usepackage{pst-node}
\usepackage{amsmath}
\usepackage{graphicx}
\usepackage{nccmath}
\journal{Journal of \LaTeX\ Templates}
\usepackage{setspace}
\usepackage{hyperref}
\usepackage{lipsum}
\usepackage{float}
\usepackage{eqnarray,amsmath}
\usepackage{amsfonts}
\usepackage{amssymb}
\usepackage{pgf}
\usepackage{graphicx}
\usepackage{grffile}
\usepackage{subfig}
\usepackage{multicol}
\usepackage{lipsum}% for dummy text
\usepackage{epstopdf}
\usepackage{mathtools}
\usepackage{mathabx}
\usepackage{pgf-pie}
\usepackage{latexsym, amsmath}
\usepackage{amssymb}
\usepackage{amsfonts}
\usepackage{amsmath}
\usepackage{amsfonts}
\usepackage{epsfig}
\numberwithin{equation}{section}
\usepackage{ragged2e}
\usepackage{pgfplots}
\usetikzlibrary{arrows, positioning, calc}
\usepackage{bchart}
\usepackage[section]{placeins}
\usepackage{subfig}
\justifying
\date{ }

\begin{document}

\begin{frontmatter}

\title{Interplay between Negation of a Probability Distribution and Jensen Inequality }

\author{Amit Srivastava}
\address{Department of Mathematics,\\ Jaypee Institute of
	Information Techno\textnormal{log}y,\\ Noida (Uttar Pradesh), India.\\
}
\begin{abstract}
Yager\cite{sp1} proposed a transformation for opposing(negating) the occurence of an event that is not certain using the idea that one can oppose the occurence of any uncertain event by allocating its probability among the other outcomes in the sample space without preference to any particular outcome \textit{i.e.} the probability of every event in the sample space is redistributed equally among the other outcomes in the sample space. However this redistribution increases the uncertainty associated with the occurence of events. In the present work, we have established bounds on the uncertainty associated with negation of a probability distribution using well known Jensen inequality. The obtained results are validated with the help of various numerical examples. Finally  a dissimilarity function between a probability distribution and its negation has been developed. 
\end{abstract}

\begin{keyword}
Negation, Uncertainty, Probability distribution, Jensen inequality, convex function.
\MSC[2010] 26B25; 94A17.
\end{keyword}
\end{frontmatter}
\section{Introduction}
Consider an $n\times n$ doubly stochastic matrix $A=(a_{ij})$ with
\[ \sum\limits_{i=1}^{n} a_{ij}=\sum\limits_{j=1}^{n} a_{ij}=1; n =2, 3, 4, ...\]
In particular, if $a_{ii} = 0 $ for all $i=1,2,\ldots,n$ and if $a_{ij} = \frac{1}{n-1} $ for all $i\neq j$; $i,j =1,2,\ldots,n$, then given a discrete finite complete probability distribution $P(n)=\{p_1,p_2,\ldots,p_{n-1},p_n\}$; its negation  $\overline{P}(n)=\{\overline{p}_1,\overline{p}_2,\ldots,\overline{p}_{n-1},\overline{p}_n\}$ \cite{sp1},\cite{sp5},\cite{sp11},\cite{sp20},\cite{sp867},\cite{spp1},\cite{spp2},\cite{spp4} can be written as
\begin{equation}
	\overline{p}_i=\sum_{j=1}^{n} a_{ij} {p}_i \;\;
\end{equation} for all $i =1,2,\ldots,n$. We can further write 
\begin{equation}
	\overline{p}_i = \frac{1}{n-1}.p_{1}+\frac{1}{n-1}.p_{2}+...+\frac{1}{n-1}.p_{i-1}+\frac{1}{n-1}.p_{i+1} +...+\frac{1}{n-1}.p_{n-1}+\frac{1}{n-1}.p_{n} = \frac{1-p_{i}}{n-1};i=1,2,\ldots,n
\end{equation}
which gives
\[\overline{p}_1 = 0.p_{1}+\frac{1}{n-1}.p_{2}+....+\frac{1}{n-1}.p_{n-1}+\frac{1}{n-1}.p_{n} = \frac{1-p_{1}}{n-1} \]
\[\overline{p}_2= \frac{1}{n-1}.p_{1}+0.p_{2}+....+\frac{1}{n-1}.p_{n-1}+\frac{1}{n-1}.p_{n} = \frac{1-p_{2}}{n-1} \]
\[.\]
\[.\]
\[.\]
\[\overline{p}_{n-1} = \frac{1}{n-1}.p_{1}+\frac{1}{n-1}.p_{2}+....+0.p_{n-1}+\frac{1}{n-1}.p_{n}= \frac{1-p_{n-1}}{n-1} \]
\[\overline{p}_n = \frac{1-p_{n}}{n-1} = \frac{1}{n-1}.p_{1}+\frac{1}{n-1}.p_{2}+\frac{1}{n-1}.p_{3}+....+0.p_{n}\]
The vector $\overline{P}(n)$ is also a probability distribution with $\overline{p}_i\in[0,1]$ \;\;for\;\; $i=1,2,\ldots,n $ and $ \sum\limits_{i=1}^{n}\overline{p}_{i}=1$. Consider the probability distribution $P(4) = \left\{\frac{1}{3} ,\frac{1}{6} ,\frac{1}{6},\frac{1}{3}\right\}$ and its corresponding negation $\overline{P}(4)= \left\{\frac{2}{9} ,\frac{5}{18} ,\frac{5}{18},\frac{2}{9}\right\}$.The first entry in $\overline{P}(4)$ is a result of equal contributions($\frac{1}{3}$ part) from the second, third and fourth entries in ${P}(4)$. Again the second entry in $\overline{P}(4)$ is a result of equal contributions($\frac{1}{3}$ part) from the first, third and fourth entries in ${P}(4)$. Similar is the case with third and fourth entries in $\overline{P}(4)$. Figure 1 shows the pictorial representation of redistribution of probabilities in $P(4)=\{p_1,p_2,p_3,p_4\}=\left\{\frac{1}{3} ,\frac{1}{6} ,\frac{1}{6},\frac{1}{3}\right\}$.If  $k$ events with zero probability are added in $P(n)$, then the revised probability distributions are 
\begin{center}
	$P(n)=\{p_{1},p_{2},...,p_{n},\underbrace{0,0,...,0}_{k-times}\}$ 
\end{center}
and
\begin{center}
	$\overline{P}(n+k) = \{\frac{1-p_{1}}{n+k-1},\frac{1-p_{2}}{n+k-1},...,\frac{1-p_{n}}{n+k-1},\underbrace{\frac{1}{n+k-1},...,\frac{1}{n+k-1}}_{k-times}\}$.  
\end{center}
 One important aspect associated with the concept of negation is that more uncertainty is inherent in it. Using the well known Shannon entropy function given by
\begin{equation}
	 H({P}(n))= H\left({p}_1,{p}_2,\ldots,{p}_n\right) = -\sum\limits _{i=1}^{n}{p}_{_{i} }{\rm log_{}\; \; }{p}_{_{i} } 
\end{equation}
we can easily check that \begin{center}
	$H({P}(4))= H(p_1,p_2,p_3,p_4)=H(\frac{1}{3} ,\frac{1}{6} ,\frac{1}{6},\frac{1}{3})<H(\frac{2}{9} ,\frac{5}{18},\frac{5}{18},\frac{2}{9})= H(\overline{p}_1,\overline{p}_2, \overline{p}_3,\overline{p}_4)=H(\overline{P}(4))$.
\end{center} This is due to the fact that it becomes more difficult to predict the occurrence of any event due to the redistribution of probabilities(In (1.3), logarithm is evaluated at base 2). Again consider the probability distributions  $P(3) = \left\{\frac{2}{3} ,\frac{1}{6} ,\frac{1}{6}\right\}$  and $Q(5) = \left\{\frac{2}{3} ,\frac{1}{6} ,\frac{1}{6},0, 0\right\}$, then its negations are given as $\overline P(3) = \left\{\frac{1}{6} ,\frac{5}{12} ,\frac{5}{12}\right\}$  and $\overline Q(5) = \left\{\frac{1}{12} ,\frac{5}{24} ,\frac{5}{24},\frac{1}{4}, \frac{1}{4}\right\}$. Here we have 
\begin{center}
$H({P}(3))= H(p_1,p_2,p_3)=H(\frac{2}{3} ,\frac{1}{6} ,\frac{1}{6})=H(\frac{2}{3} ,\frac{1}{6},\frac{1}{6},0, 0)= H(q_1,q_2,q_3,q_4,q_5)=H({Q}(5))$.
\end{center}
but  
\begin{center}
	$H(\overline{P}(3))= H(\overline{p}_1,\overline{p}_2, \overline{p}_3)=H(\frac{1}{6} ,\frac{5}{12} ,\frac{5}{12})<H(\frac{1}{12} ,\frac{5}{24},\frac{5}{24},\frac{1}{4},\frac{1}{4})= H(\overline{q}_1,\overline{q}_2, \overline{q}_3,\overline{q}_4,\overline{q}_5)=H(\overline{Q}(5))$.
\end{center}
The Shannon uncertainty associated with ${P}(3)$ and ${Q}(5)$ is identical since an event with zero probability will not contribute anything to the total uncertainty embedded in ${Q}(5)$, but the uncertainty associated with $(\overline{Q}(5))$ is greater than that of $(\overline{P}(3))$ since the Yager's negation defined by (1.1) gives identical weightage to events with zero and non zero probabilities. Any probability distribution ${P}(n)$ and its negation $\overline{P}(n)$ will be identical only when all the probabilities in ${P}(n)$ are equal since in that case any further redistrbution will not alter the entries of ${P}(n)$. Another well known fact about the Shannon function given by (1.3) is that
\begin{equation}
H({P}(n))= H\left({p}_1,{p}_2,\ldots,{p}_n\right) \leq {\rm log_{}\;}n 
\end{equation} 
where $n \in \mathbb{N}$. The inequality(1.4) can be easily established using the well known Shannon inequality given as
\begin{equation}
-\sum_{i=1}^{n}p_i {\rm log_{}\;}  p_i \leq -\sum_{i=1}^{n}p_i {\rm log_{}\;} q_i
\end{equation}
Here $0 \leq p_i\leq 1; 0 \leq q_i\leq 1; 1 \leq i\leq n  $ and $\sum_{i=1}^{n}p_i =1= \sum_{i=1}^{n}q_i$.
Substituting $q_i =\frac{1}{n}; 1 \leq i\leq n$ in (1.5) yields (1.4). Again replacing $p_i$ by $\overline p_i$ and substituting $q_i =\frac{1}{n}; 1 \leq i\leq n$ in (1.5) gives 
\begin{equation}
	H(\overline{P}(n))= H\left(\overline{p}_1,\overline{p}_2,\ldots,\overline{p}_n\right) \leq {\rm log_{}\;}n 
\end{equation} 
Since $\overline{P}(n)$ is more uncertain than ${P}(n)$, therefore from (1.4) and (1.6), we can write
\begin{equation}
{\rm log_{}\;}n - H(\overline{P}(n))\leq {\rm log_{}\;}n - H({P}(n))
\end{equation}
Equality holds in (1.4), (1.5), (1.6) and (1.7) when all the probabilties in $P(n)=\{p_1,p_2,\ldots,p_n\}$ are equal. Yager\cite{sp1} also defined the general form of double negation as 
\[\bar{\bar{p}}_{i} = \frac{1}{n-1}.\overline p_{1}+\frac{1}{n-1}.\overline p_{2}+\frac{1}{n-1}.\overline p_{3}+...+\frac{1}{n-1}.\overline p_{i-1}+\frac{1}{n-1}.\overline p_{i+1} +...+\frac{1}{n-1}.\overline p_{n}\]
\begin{equation}
	= \frac{1-\overline{p}_{i} }{n-1} =\frac{1-\left(\frac{1-p_{i} }{n-1} \right)}{n-1} =\frac{p_{i} +n-2}{\left(n-1\right)^{2} }; i=1,2,\ldots,n\
\end{equation}
It is clear that Yager's double negation is a result of redistribution of entries in $\overline{P}(n)=\{\overline{p}_1,\overline{p}_2,\ldots,\overline{p}_n\}$. Again considering $P(4) = \left\{\frac{1}{3} ,\frac{1}{6} ,\frac{1}{6},\frac{1}{3}\right\}$; $\overline{P}(4)= \left\{\frac{2}{9} ,\frac{5}{18} ,\frac{5}{18},\frac{2}{9}\right\}$  and $\overline{\overline{P}}(4)= \left\{\frac{7}{27} ,\frac{13}{54} ,\frac{13}{54},\frac{7}{27}\right\}$ and observing that 
\begin{center}
	$H({P}(4))= H(p_1,p_2,p_3,p_4)=H(\frac{1}{3} ,\frac{1}{6} ,\frac{1}{6},\frac{1}{3})<H(\frac{2}{9} ,\frac{5}{18},\frac{5}{18},\frac{2}{9})= H(\overline{p}_1,\overline{p}_2, \overline{p}_3,\overline{p}_4)=H(\overline{P}(4))$
\end{center}
\begin{center}
	$<H(\frac{7}{27} ,\frac{13}{54},\frac{13}{54},\frac{13}{54})= H(\overline{\overline{p}}_1,\overline{\overline{p}}_2, \overline{\overline{p}}_3,\overline{\overline{p}}_4)=H(\overline{\overline{P}}(4))$.
\end{center}
we conclude that 
\begin{equation}
	H({P}(n))\leq H(\overline{P}(n)) \leq H(\overline{\overline{P}}(n))\leq {\rm log_{}\;}n 
\end{equation}
Some generalization and extensions of Yager's negation can be seen in  \cite{sp14},\cite{sp66},\cite{sp15},\cite{sp26}, \cite{sp36},\cite{sp46} and \cite{sp56}.
It is now clear that negation of a probability distribution opposes the information inherent in the distribution. In fact, if a probability distribution is uncertain(a state other than maximum entropy), the more iterations of negation, the more uncertain this probability event becomes eventually converging to a homogeneous state\textit{ i.e.} maximum entropy. In this respect, it is identical to many natural phenomena which tend towards maximum entropy arrangements with passage of time. Therefore investigating the basic structure of negation in probabilistic frameworks may result in new insights as far as occurrence of natural phenomena is concerned. In section 2, we have established some inequalities involving $P(n)$ and $\overline{P}(n)$ which are shown to be important from the information theoretic point of view. In Section 3, a new dissimilarity function between a probability distribution and its negation has been developed. Section 4 concludes this paper.
\section{Main Result}
We begin by stating the Jensen inequality (without proof). Consider a convex function $f: I \rightarrow \mathbb{R} $, where $I$ is an interval in $\mathbb{R}$ and $x_{1},x_{2},...,x_{n}$ are in $I$. If $\alpha_{1},\alpha_{2},...,\alpha_{n}$ are non-negative real numbers such that $ \sum_{1}^{n} \alpha_{i} =1 $, then we have
\begin{equation}
	g(\alpha_{1}x_{1} + \alpha_{2}x_{2} + ... + \alpha_{n}x_{n}) \leq \alpha_{1}g(x_{1}) +\alpha_{2}g(x_{2})+...\\+\alpha_{n}g(x_{n})
\end{equation}
If $f$ is concave then the above inequality will be reversed. Equality holds in (2.1) if and only if all  entries in the set $\left(x_{1},x_{2},...,x_{n}\right)$ are equal. Now we will utilize the above inequality for determining various information theoretic inequalties involving $P(n)$ and $\overline{P}(n)$.\\ 
\begin{figure}[ht!]
	\centering
	\subfloat{\resizebox{0.4\textwidth}{!}{\begin{tikzpicture}[every node/.style={circle,draw=black,align=center}] \node(c){$0 +\frac{1}{18}+\frac{1}{18}+\frac{1}{9} = \frac{2}{9}=\overline{p}_1$};\node(a)[above left =of c]{$0.p_1 = 0.\frac{1}{3}$} edge [->](c); \node(d)[left= of c]{$\frac{1}{3}.p_2 =\frac{1}{3}.\frac{1}{6} $} edge [->] (c);\node(b)[below left =of c]{$\frac{1}{3}.p_3 =\frac{1}{3}.\frac{1}{6} $} edge [->](c);\node(b)[below =of b]{$\frac{1}{3}.p_4 =\frac{1}{3}.\frac{1}{3} $} edge [->](c) ;\end{tikzpicture}}}
	\quad
	\subfloat{\resizebox{0.4\textwidth}{!}{\begin{tikzpicture}[every node/.style={circle,draw=black,align=center}] \node(c){$\frac{1}{9} + 0 +\frac{1}{18}+\frac{1}{9} = \frac{5}{18}=\overline{p}_2$};\node(a)[above left =of c]{$\frac{1}{3}.p_1 = \frac{1}{3}.\frac{1}{3}$} edge [->](c); \node(d)[left= of c]{$0.p_2 =0.\frac{1}{6} $} edge [->] (c);\node(b)[below left =of c]{$\frac{1}{3}.p_3 =\frac{1}{3}.\frac{1}{6} $} edge [->](c);\node(b)[below =of b]{$\frac{1}{3}.p_4 =\frac{1}{3}.\frac{1}{3} $} edge [->](c) ;\end{tikzpicture}}} \\
	\subfloat{\resizebox{0.4\textwidth}{!}{\begin{tikzpicture}[every node/.style={circle,draw=black,align=center}] \node(c){$\frac{1}{9} +\frac{1}{18} + 0 +\frac{1}{9} = \frac{5}{18}=\overline{p}_3$};\node(a)[above left =of c]{$\frac{1}{3}.p_1 = \frac{1}{3}.\frac{1}{3}$} edge [->](c); \node(d)[left= of c]{$\frac{1}{3}.p_2 =\frac{1}{3}.\frac{1}{6} $} edge [->] (c);\node(b)[below left =of c]{$0.p_3 =0.\frac{1}{6} $} edge [->](c);\node(b)[below =of b]{$\frac{1}{3}.p_4 =\frac{1}{3}.\frac{1}{3} $} edge [->](c) ;\end{tikzpicture}}} \quad
	\subfloat{\resizebox{0.4\textwidth}{!}{\begin{tikzpicture}[every node/.style={circle,draw=black,align=center}] \node(c){$\frac{1}{9} +\frac{1}{18} +\frac{1}{18}+0 = \frac{2}{9}=\overline{p}_4$};\node(a)[above left =of c]{$\frac{1}{3}.p_1 = \frac{1}{3}.\frac{1}{3}$} edge [->](c); \node(d)[left= of c]{$\frac{1}{3}.p_2 =\frac{1}{3}.\frac{1}{6} $} edge [->] (c);\node(b)[below left =of c]{$\frac{1}{3}.p_3 =\frac{1}{3}.\frac{1}{6} $} edge [->](c);\node(b)[below =of b]{$0.p_4 =0.\frac{1}{2} $} edge [->](c) ;\end{tikzpicture}} }
	\caption{Negation of $P(4)=\{p_1,p_2,p_3,p_4\}=\left\{\frac{1}{3} ,\frac{1}{6} ,\frac{1}{6},\frac{1}{3}\right\}$ }
\end{figure}
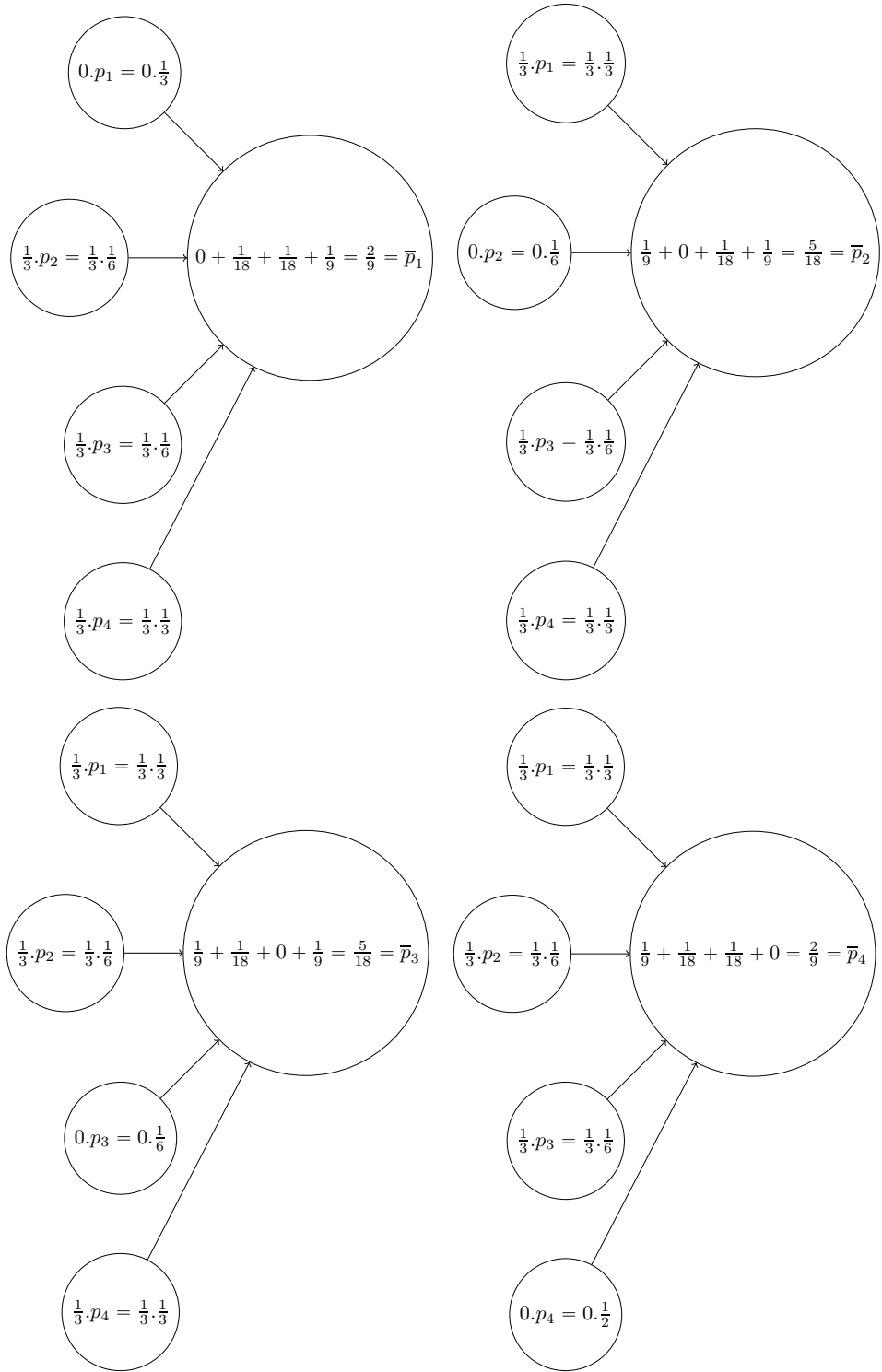

\textbf{Theorem 2.1}: Consider a convex function $f: I \rightarrow \mathbb{R} $, where $I$ is an interval in $\mathbb{R}$. For the probability distribution $P(n)=\{p_1,p_2,\ldots,p_n\}$ and its negation $\overline{P}(n)=\{\overline{p}_1,\overline{p}_2,\ldots,\overline{p}_n\}$, we have 
\begin{equation}
	f\left(\frac{1}{n}\right)\leq \left(\frac{1}{n^2}\right)   \sum\limits_{i=1}^{n}f({p}_{i})+\left(\frac{n-1}{n^2}\right)\sum\limits_{i=1}^{n}f(\overline{p}_{i})
\end{equation}
Equality holds in (2.2) when all the probabilties in $P(n)=\{p_1,p_2,\ldots,p_n\}$ are equal.\\
\textbf{Proof}: For $i = 1,2,\ldots,n$, we have
\[f\left(\frac{1}{n}\right)=f\left(\frac{{{p}_1 + {p}_2 + ...+{p}_{i-1}+{p}_{i}+{p}_{i+1} +...+ {p}_n}}{n}\right)\]
\[=f\left(\frac{{{p}_i}}{n}+\left(\frac{n-1}{n}\right)\left(\frac{{p}_1+ {p}_2 + ...+{p}_{i-1}+{p}_{i+1} + ...+{p}_n}{n-1}\right)\right)\]
\[=f\left(\frac{{{p}_i}}{n}+\left(\frac{n-1}{n}\right)\overline{p}_i\right)\]  
\begin{equation}
\leq \frac{1}{n} f\left({p}_i\right)+\left(\frac{n-1}{n}\right) f\left(\overline{p}_i\right)
\end{equation}
The above follows directly from (2.1). Summing over all $i$ in (2.3) gives (2.2). Taking $f(x) = -{\rm log_{}\;} x = \Delta(x)$(say) in (2.3), we obtain
\begin{equation}
	{\rm log_{}\;} n \leq \frac{1}{n} \left(-{\rm log_{}\;} {p}_i\right) +\left(\frac{n-1}{n}\right) \left(-{\rm log_{}\;} \overline {p}_i\right)
\end{equation} 
for all $i = 1,2,\ldots,n$. Here $\Delta({p}_i) = -{\rm log_{}\;} {p}_i$ denotes the self information associated with an event occuring with probability ${p}_i(0\leq{p}_i\leq1)$. The inequality (2.4) shows that maximum value of the Shannon entropy function given by (1.3) is less than than the convex combination of self information associated with events having probabilities ${p}_i(0\leq{p}_i\leq1)$ and  $\overline{p}_i(0\leq \overline{p}_i\leq \frac{1}{n-1})$;$i = 1,2,\ldots,n$. Taking $f(x) = -{\rm log_{}\;} x = \Delta(x)$(say) in (2.2), we obtain 
\begin{equation}
	{\rm log_{}\;} n \leq \frac{1}{n^2} \sum\limits_{i=1}^{n} \left(\Delta {p}_i\right) +\left(\frac{n-1}{n^2}\right) \sum\limits_{i=1}^{n}\left(\Delta \overline {p}_i\right)
\end{equation}
\textbf{Theorem 2.2}: Consider a convex function $f: I \rightarrow \mathbb{R} $, where $I$ is an interval in $\mathbb{R}$. For the probability distribution $\overline{P}(n)=\{\overline{p}_1,\overline{p}_2,\ldots,\overline{p}_n\}$ and its negation $\overline{\overline{P}}(n)=\{{\overline{\overline{p}}_{1},\overline{\overline{p}}_{2},...,\overline{\overline{p}}_{n}}\}$, we have 
\begin{equation}
	f\left(\frac{1}{n}\right)\leq \left(\frac{1}{n^2}\right)   \sum\limits_{i=1}^{n}f(\overline{p}_{i})+\left(\frac{n-1}{n^2}\right)\sum\limits_{i=1}^{n}f(\overline{\overline{p}}_{i})
\end{equation}
Equality holds in (2.6) when all the probabilties in $P(n)=\{p_1,p_2,\ldots,p_n\}$ are equal.\\
\textbf{Proof}: Using (1.8), the proof follows on similar lines as the proof of theorem (2.1).\\
\textbf{Theorem 2.3}: Consider a concave function $f: I \rightarrow \mathbb{R} $, where $I$ is an interval in $\mathbb{R}$. For the probability distribution $P(n)=\{p_1,p_2,\ldots,p_n\}$ and its negation $\overline{P}(n)=\{\overline{p}_1,\overline{p}_2,\ldots,\overline{p}_n\}$, we have 
\begin{equation}
	\left(\frac{1}{n^2}\right)   \sum\limits_{i=1}^{n}f({p}_{i})+\left(\frac{n-1}{n^2}\right)\sum\limits_{i=1}^{n}f(\overline{p}_{i})\leq f\left(\frac{1}{n}\right)
\end{equation}
Equality holds in (2.7) when all the probabilties in $P(n)=\{p_1,p_2,\ldots,p_n\}$ are equal.\\
\textbf{Proof}:Using inequality (2.1) for concave functions, the proof follows on similar lines as the proof of theorem (2.1). Taking $f(x) = -x{\rm log_{}\;} x $ in (2.7), we obtain
\begin{center}
	$ -\left(\frac{1}{n^2}\right)\sum\limits_{i=1}^{n}{p}_{i} {\rm log_{}}{p}_{i}-\left(\frac{n-1}{n^2}\right)\sum\limits_{i=1}^{n}\overline{p}_{i} {\rm log_{}}\overline{p}_{i} \leq -\left(\frac{1}{n}\right) {\rm log_{}}\left(\frac{1}{n}\right)$
\end{center}
which gives
\begin{center}
	$\left(\frac{1}{n}\right)H({P}(n))+\left(\frac{n-1}{n}\right)H(\overline{P}(n))\leq{\rm log_{}} n$
\end{center}
which is in agreement with (1.4) and (1.6). However for n=2, the above reduces to (1.4) or (1.6) since for any probability distribution $P(2)=\{p_1,p_2\}$, its negation will be $\overline P(2)=\{p_2,p_1\}$.\\
\textbf{Theorem 2.4}: Consider a convex $f: I \rightarrow \mathbb{R} $, where $I$ is an interval in $\mathbb{R}$ and the probability distributions $P(n)=\{p_1,p_2,\ldots,p_{n-2},p_{n-1},p_n\}$ and  $\overline{P}(n)=\{\overline{p}_1,\overline{p}_2,\ldots,\overline{p}_{n-2},\overline{p}_{n-1},\overline{p}_n\}$. For $i = 1,2,\ldots,n$, we have
\[f\left(\zeta_{n-1}\left({{p}_1,{p}_2,...,{p}_{i-1},{p}_{i+1},...,p_{n-2},p_{n-1},{p}_n}\right)\right)\]
\[\leq \left(\frac{1}{n-1}\right)f\left({p}_{n}\right)+\left(\frac{n-2}{n-1}\right)f\left(\zeta_{n-2}\left({{p}_1,{p}_2,...,{p}_{i-1},{p}_{i+1},...,p_{n-2},p_{n-1}}\right)\right)\]
\[\leq \left(\frac{1}{n-1}\right)f\left({p}_{n}\right)+\left(\frac{n-2}{n-1}\right)\left(\frac{1}{n-1}\right)f\left({p}_{n-1}\right)\]
\[+ \left(\frac{n-2}{n-1}\right)\left(\frac{n-3}{n-2}\right)f\left(\zeta_{n-3}\left({{p}_1,{p}_2,...,{p}_{i-1},{p}_{i+1},...,p_{n-2}}\right)\right)\]
\[...\]
\[\leq \left(\frac{1}{n-1}\right)f\left({p}_{n}\right) + \left(\frac{1}{n-1}\right)f\left({p}_{n-1}\right)+\left(\frac{1}{n-1}\right)f\left({p}_{n-2}\right)+...+\left(\frac{1}{n-1}\right)f\left({p}_{i-1}\right)\]
\begin{equation}
+\left(\frac{1}{n-1}\right)f\left({p}_{i+1}\right)...+\left(\frac{1}{n-1}\right)f\left({p}_{2}\right)+\left(\frac{1}{n-1}\right)f\left({p}_{1}\right)\
\end{equation}
Here 
$\zeta_{n-1}\left({{p}_1,{p}_2,{p}_{i-1},...,{p}_{i+1},...,p_{n-2},p_{n-1},{p}_n}\right)$
\[=\frac{{{p}_1 + {p}_2 + ...+{p}_{i-1}+{p}_{i+1} +...+p_{n-2}+p_{n-1}+{p}_n}}{n-1}\];
$\zeta_{n-2}\left({{p}_1,{p}_2,{p}_{i-1},...,{p}_{i+1},...,p_{n-2},p_{n-1}}\right)$
\[=\frac{{{p}_1 + {p}_2 + ...+{p}_{i-1}+{p}_{i+1} +...+p_{n-2}+p_{n-1}}}{n-2}\]
and so on. Equality holds in (2.8) when all the probabilties in $P(n)=\{p_1,p_2,\ldots,p_n\}$ are equal. \\
\textbf{Proof}:For $i = 1,2,\ldots,n$, we have
\[f\left(\overline{p}_{i}\right)=f\left(\zeta_{n-1}\left({{p}_1,{p}_2,...,{p}_{i-1},{p}_{i+1},...,p_{n-2},p_{n-1},{p}_n}\right)\right)\]
\[=f\left(\frac{{{p}_1 + {p}_2 + ...+{p}_{i-1}+{p}_{i+1} +...+p_{n-2}+p_{n-1}+{p}_n}}{n-1}\right)\]
\[=f\left(\frac{{{p}_n}}{n-1}+\left(\frac{n-2}{n-1}\right)\left(\frac{{{p}_1 + {p}_2 + ...+{p}_{i-1}+{p}_{i+1} +...+,p_{n-2}+p_{n-1}}}{n-2}\right)\right)\]
\[=f\left(\frac{{{p}_n}}{n-1}+\left(\frac{n-2}{n-1}\right)\zeta_{n-2}\left({{p}_1,{p}_2,{p}_{i-1},...,{p}_{i+1},...,p_{n-2},p_{n-1}}\right)\right)\]
\[\leq \left(\frac{1}{n-1}\right)f\left({p}_{n}\right)+\left(\frac{n-2}{n-1}\right)f\left(\zeta_{n-2}\left({{p}_1,{p}_2,...,{p}_{i-1},{p}_{i+1},...,p_{n-2},p_{n-1}}\right)\right)\]
Again
\[f\left(\zeta_{n-2}\left({{p}_1,{p}_2,...,{p}_{i-1},{p}_{i+1},...,p_{n-2},p_{n-1}}\right)\right)\]
\[=f\left(\frac{{{p}_1 + {p}_2 + ...+{p}_{i-1}+{p}_{i+1} +...+p_{n-2}+p_{n-1}}}{n-2}\right)\]
\[=f\left(\frac{{p_{n-1}}}{n-2}+\left(\frac{n-3}{n-2}\right)\left(\frac{{{p}_1 + {p}_2 + ...+{p}_{i-1}+{p}_{i+1} +...+p_{n-2}}}{n-3}\right)\right)\]
\[=f\left(\frac{{p_{n-1}}}{n-2}+\left(\frac{n-3}{n-2}\right)\zeta_{n-3}\left({{p}_1,{p}_2,...,{p}_{i-1},{p}_{i+1},...,p_{n-2}}\right)\right)\]
\[\leq \left(\frac{1}{n-2}\right)f\left({p}_{n-1}\right)+\left(\frac{n-3}{n-2}\right)f\left(\zeta_{n-3}\left({{p}_1,{p}_2,...,{p}_{i-1},{p}_{i+1},...,p_{n-2}}\right)\right)\]
Therefore 
\[f\left(\overline{p}_{i}\right)\leq \left(\frac{1}{n-1}\right)f\left({p}_{n}\right)+\left(\frac{n-2}{n-1}\right)\left(\frac{1}{n-1}\right)f\left({p}_{n-1}\right)\]
\[+ \left(\frac{n-2}{n-1}\right)\left(\frac{n-3}{n-2}\right)f\left(\zeta_{n-3}\left({{p}_1,{p}_2,...,{p}_{i-1},{p}_{i+1},...,p_{n-2}}\right)\right)\]
\[= \left(\frac{1}{n-1}\right)f\left({p}_{n}\right)+\left(\frac{1}{n-1}\right)f\left({p}_{n-1}\right)\]
\[+ \left(\frac{n-3}{n-1}\right)f\left(\frac{{{p}_1 + {p}_2 + ...+{p}_{i-1}+{p}_{i+1} +...+p_{n-2}}}{n-3}\right)\]
\[...\]
\[\leq \left(\frac{1}{n-1}\right)f\left({p}_{n}\right) + \left(\frac{1}{n-1}\right)f\left({p}_{n-1}\right)+\left(\frac{1}{n-1}\right)f\left({p}_{n-2}\right)+...+\left(\frac{1}{n-1}\right)f\left({p}_{i-1}\right)\]
\[+\left(\frac{1}{n-1}\right)f\left({p}_{i+1}\right)...+\left(\frac{1}{n-1}\right)f\left({p}_{2}\right)+\left(\frac{1}{n-1}\right)f\left({p}_{1}\right)\]
which completes the proof. Using (2.8), for $i = 1,2,\ldots,n$, we can write 
\begin{equation}
	f\left(\overline{p}_{i}\right)\leq \left(\frac{1}{n-1}\right)\sum\limits_{j=1, i\neq j}^{n}f\left({p}_{j}\right)
\end{equation}
Taking $f(x) = -{\rm log_{}\;} x = \Delta(x)$(say) in (2.9), we obtain
\begin{equation}
	\left(-{\rm log_{}\;} \overline {p}_i\right)\leq \left(\frac{1}{n-1}\right)\sum\limits_{j=1, i\neq j}^{n}\left(-{\rm log_{}\;}{p}_i\right)
\end{equation}
Equality holds in above when all the probabilties in $P(n)=\{p_1,p_2,\ldots,p_n\}$ are equal. However we consider a specific case here. Consider a random variable $X=(x_{1},x_{2},x_{3}, x_{4},x_{5})$ with corresponding probabilities $P(5) = (p(x_{1}),p(x_{2}),p(x_{3}),p(x_{4}),p(x_{5}))$. We assume that the probability distribution $P(5)$ 
is symmetric about $X = x_{3}$ and the probabilities in $P(5)$ satisfy
\[p(x_{1})+ p(x_{2})= \frac{p(x_{3})}{2}=p(x_{4})+ p(x_{5})\]
\[p(x_{1})= p(x_{2}); p(x_{4})= p(x_{5})\]
which gives 
\[p(x_{1})= \frac{1}{8}, p(x_{2})= \frac{1}{8},p(x_{3})= \frac{1}{2},p(x_{4})= \frac{1}{8},p(x_{5})= \frac{1}{8}.\]
Now the negation of $p(x_{3})$ can be written as 
\[\overline p(x_{3})= \frac{1-p(x_{3})}{4}=\frac{p(x_{1}) + p(x_{2}) + p(x_{4}) + p(x_{5})}{4}=\frac{p(x_{1}) + p(x_{2})}{2}\]
Therefore we finally have 
\[\overline p(x_{3})=p(x_{1})=p(x_{2})=p(x_{4})=p(x_{5})\]
Now for the probability distribution $P(5)$, writing(2.10) for $i=3$ gives
\begin{center}
	$\left(-{\rm log_{}\;} \overline {p}_3\right)= \left(\frac{1}{4}\right)\left(-{\rm log_{}\;} {p}_1-{\rm log_{}\;} {p}_2-{\rm log_{}\;} {p}_4-{\rm log_{}\;} {p}_5\right)$
\end{center}
Therefore equality holds in (2.10) for a probability distribution that is symmetric but not uniform. However the equality holds about the point where the distribution is symmetric.  The same result can be verified for (2.8) also. It is clear from Figure 2 that negation transformation has preserved the symmetry of $P(5)$. 
\begin{figure}
	\centering
	\resizebox{0.3\textwidth}{!}{\begin{tikzpicture}
5			\pie[color={black!12.5, black!12.5, black!50, black!12.5, black!12.5}]{12.5/${p}_1$, 12.5/${p}_2$, 50/${p}_3$, 12.5/${p}_4$, 12.5/${p}_5$}
	\end{tikzpicture}}
	\resizebox{0.3\textwidth}{!}{
		\begin{tikzpicture}
			\pie[color={black!21.875, black!21.875,black!12.5, black!21.875, black!21.875}]
			{21.875/$\overline{p}_1$, 21.875/$\overline{p}_2$, 12.5/ $\overline{p}_3$, 21.875/$\overline{p}_4$, 21.875/$\overline{p}_5$}		
	\end{tikzpicture}}
	\resizebox{0.3\textwidth}{!}{\begin{tikzpicture}
			\pie[color={black!20, black!20,black!20,black!20, black!20}]
			{20/${p}_1^*$, 20/${p}_2^*$,20/ ${p}_3^*$, 20/${p}_4^*$, 20/${p}_5^*$}		
	\end{tikzpicture}}
	\caption{Probability distribution $P(5) = \left\{ \frac{1}{8},  \frac{1}{8}, \frac{1}{2},  \frac{1}{8},  \frac{1}{8}\right\}$,$\overline {P}(5) = (0.225, 0.2,
		0.15, 0.2, 0.225)$ and the MEPD $P^*(5)= \left\{\frac{7}{32}, \frac{7}{32}, \frac{1}{8}, \frac{7}{32}, \frac{7}{32}\right\}$ }
\end{figure}
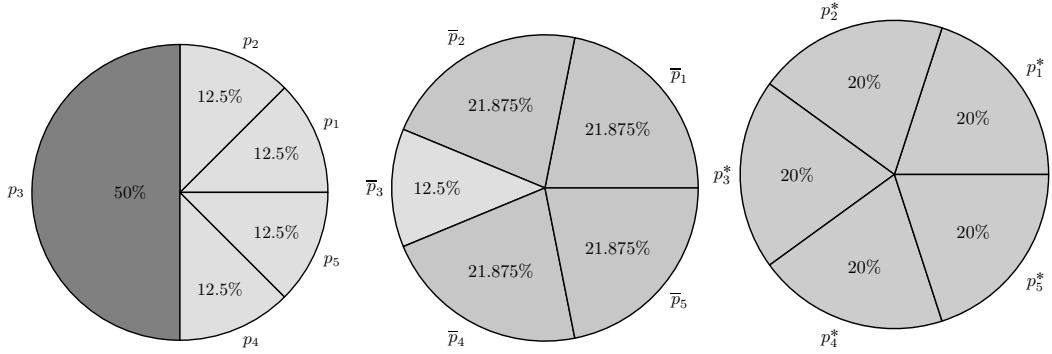
\section{Quantifying the dissimilarity between a probability distribution and its negation}
From (1.7), it is clear dissimilarity between a probability distribution and its negation signifies the dissimilarity between the probability distribution and the uniform distribution which is the maximum  entropy probability distribution in the absence of any additional information regarding the occurence of events. So quantifying dissimilarity between the probability distribution is important from the information theoretic point of view \cite{spp5}, \cite{spp6}. We start with following theorem. \\
\textbf{Theorem 3.1.} For all  $P(n)=\{p_1,p_2,\ldots,p_n\}$ and $\overline{P}(n)=\{\overline{p}_1,\overline{p}_2,\ldots,\overline{p}_n\}$, the function defined by 
\begin{equation}
	I_{\alpha } (P(n);\overline{P}(n))=-log\left(\frac{1+\frac{1}{2} \sum _{i=1}^{n}\left(min\left(p_{i} ,\frac{\left(2^{\alpha } -1\right)p_{i} +q_{i} }{2^{\alpha } } \right)+min\left(\frac{p_{i} +\left(2^{\alpha } -1\right)q_{i} }{2^{\alpha } } ,q_{i} \right)\right) }{2} \right)\, \, \, ;\alpha =0,1,2,3,...
\end{equation}
represents a series of dissimilarity functions between $P(n)$ and  $\overline{P}(n)$ satisfying the following properties

\noindent \textbf{M1. } $0\le I_{\alpha } (P(n);\overline{P}(n))\le 1$ . Further${\ I}_{\alpha }\left(P(n);\overline{P}(n)\right)=0$ if and only if $P(n)=\overline{P}(n)$.

\noindent \textbf{M2.     }$I_{\alpha } (P(n);\overline{P}(n))=\, \, \, I_{\alpha } (\overline{P}(n);P(n))$

\noindent \textbf{M3.     }$I_{0} ((P(n);\overline{P}(n))\le I_{1} (P(n);\overline{P}(n))\le I_{2}(P(n);\overline{P}(n))\le ...\le I_{\alpha } (P(n);\overline{P}(n))\le ...\forall i$\textbf{  }

\noindent \textbf{}

\noindent \textbf{Proof.  }For all  $P(n)=\{p_1,p_2,\ldots,p_n\}$ and $\overline{P}(n)=\{\overline{p}_1,\overline{p}_2,\ldots,\overline{p}_n\}$, we have
\noindent 
\[min\left(p_{i} ,\overline{p}_{i} \right)\le \left(\frac{p_{i} +\overline{p}_{i} }{2} \right)\forall i                   \Rightarrow 0\le \sum _{i=1}^{n}min\left(p_{i} ,\overline{p}_{i} \right) \le \sum _{i=1}^{n}\left(\frac{p_{i} +\overline{p}_{i} }{2} \right) =1\] 
\[\Rightarrow \frac{1}{2} \le \left(\frac{1+\sum _{i=1}^{n}min\left(p_{i} ,\overline{p}_{i} \right) }{2} \right)\le 1            \Rightarrow 0\le -log\left(\frac{1+\sum _{i=1}^{n}min\left(p_{i} ,\overline{p}_{i} \right) }{2} \right)\le 1\] 
\[\Rightarrow 0\le I_{0} (P(n);\overline{P}(n))\le 1\] 

\noindent Again for all $P(n)=\{p_1,p_2,\ldots,p_n\}$ and $\overline{P}(n)=\{\overline{p}_1,\overline{p}_2,\ldots,\overline{p}_n\}$, we have

\noindent         $min\left(p_{i} ,\left(\frac{p_{i} +\overline{p}_{i} }{2} \right)\right)\le \frac{p_{i} +\left(\frac{p_{i} +\overline{p}_{i} }{2} \right)}{2} $ and  $min\left(\left(\frac{p_{i} +\overline{p}_{i} }{2} \right),\overline{p}_{i} \right)\le \frac{\left(\frac{p_{i} +\overline{p}_{i} }{2} \right)+\overline{p}_{i} }{2} \forall i$

\noindent $\Rightarrow 0\le \sum _{i=1}^{n}min\left(p_{i} ,\left(\frac{p_{i} +\overline{p}_{i} }{2} \right)\right) \le \sum _{i=1}^{n}\left(\frac{p_{i} +\left(\frac{p_{i} +\overline{p}_{i} }{2} \right)}{2} \right) =1$ and
\[0\le \sum _{i=1}^{n}min\left(\left(\frac{p_{i} +\overline{p}_{i} }{2} \right),\overline{p}_{i} \right) \le \sum _{i=1}^{n}\left(\frac{\left(\frac{p_{i} +\overline{p}_{i} }{2} \right)+\overline{p}_{i} }{2} \right) =1\] 
\[\Rightarrow 0\le \frac{1}{2} \left(\sum _{i=1}^{n}min\left(p_{i} ,\left(\frac{p_{i} +\overline{p}_{i} }{2} \right)\right) +\sum _{i=1}^{n}min\left(\left(\frac{p_{i} +\overline{p}_{i} }{2} \right),\overline{p}_{i} \right) \right)\le 1\] 
\[\Rightarrow 0\le -log\left(\frac{1+\frac{1}{2} \sum _{i=1}^{n}\left(min\left(p_{i} ,\left(\frac{p_{i} +\overline{p}_{i} }{2} \right)\right)+min\left(\left(\frac{p_{i} +\overline{p}_{i} }{2} \right),\overline{p}_{i} \right)\right) }{2} \right)\le 1\] 
\[\Rightarrow 0\le I_{1} (P(n);\overline{P}(n))\le 1\] 

\noindent Again, a choice of pairs $\left({p}_{i}, \frac{3p_{i} +q_{i} }{4} \right)$  and $\left({\frac{p_i+3\overline{p}_{i}}{4},\overline{p}_{i}}\right)$ will give 

\noindent 
\[\Rightarrow 0\le -log\left(\frac{1+\frac{1}{2} \sum _{i=1}^{n}\left(min\left(p_{i} ,\left(\frac{3p_{i} +\overline{p}_{i} }{4} \right)\right)+min\left(\left(\frac{p_{i} +3\overline{p}_{i} }{4} \right),\overline{p}_{i} \right)\right) }{2} \right)\le 1\] 
\[\Rightarrow 0\le I_{2} (P(n);\overline{P}(n))\le 1\]

\noindent Further a choice of pairs  $\left(p_i,\frac{{7p}_i+\overline{p}_{i}}{8}\right)$ and $\left({\frac{p_i+7\overline{p}_{i}}{8},\overline{p}_{i}}\right)$ will give 

\noindent 
\[\Rightarrow 0\le -log\left(\frac{1+\frac{1}{2} \sum _{i=1}^{n}\left(min\left(p_{i} ,\left(\frac{7p_{i} +\overline{p}_{i} }{8} \right)\right)+min\left(\left(\frac{p_{i} +7\overline{p}_{i} }{8} \right),\overline{p}_{i} \right)\right) }{2} \right)\le 1\] 
\[\Rightarrow 0\le I_{3} (P(n);\overline{P}(n))\le 1\] 
Finally a choice of pairs  $\left(p_i,\frac{{\left(2^{\propto }-1\right)p}_i+\overline{p}_{i}}{2^{\alpha }}\right)$ and $\left({\frac{p_i+\left(2^{\propto }-1\right)\overline{p}_{i}}{2^{\alpha }},\overline{p}_{i}}\right)$ will give 

\noindent 
\[\Rightarrow 0\le -log\left(\frac{1+\frac{1}{2} \sum _{i=1}^{n}\left(min\left(p_{i} ,\left(\frac{\left(2^{\alpha } -1\right)p_{i} +\overline{p}_{i} }{2^{\alpha } } \right)\right)+min\left(\left(\frac{p_{i} +\left(2^{\alpha } -1\right)\overline{p}_{i} }{2^{\alpha } } \right),\overline{p}_{i} \right)\right) }{2} \right)\le 1\] 
\[\Rightarrow 0\le I_{\alpha } (P(n);\overline{P}(n))\le 1\, \, \, ;\alpha =0,1,2,3,...\] 
Further${\ I}_{\alpha }\left(P(n);\overline{P}(n)\right)=0$ implies\textbf{ }
\[-log\left(\frac{1+\frac{1}{2} \sum _{i=1}^{n}\left(min\left(p_{i} ,\left(\frac{\left(2^{\alpha } -1\right)p_{i} +\overline{p}_{i} }{2^{\alpha } } \right)\right)+min\left(\left(\frac{p_{i} +\left(2^{\alpha } -1\right)\overline{p}_{i} }{2^{\alpha } } \right),\overline{p}_{i} \right)\right) }{2} \right)=0=-log1\] 
which implies
\[\sum _{i=1}^{n}\left(min\left(p_{i} ,\left(\frac{\left(2^{\alpha } -1\right)p_{i} +\overline{p}_{i} }{2^{\alpha } } \right)\right)+min\left(\left(\frac{p_{i} +\left(2^{\alpha } -1\right)\overline{p}_{i} }{2^{\alpha } } \right),\overline{p}_{i} \right)\right) =2\] 
\[\Rightarrow \frac{1}{2} \sum _{i=1}^{n}\left(\begin{array}{l} {\left(p_{i} +\left(\frac{\left(2^{\alpha } -1\right)p_{i} +\overline{p}_{i} }{2^{\alpha } } \right)\right)-\left|p_{i} -\left(\frac{\left(2^{\alpha } -1\right)p_{i} +\overline{p}_{i} }{2^{\alpha } } \right)\right|} \\ {+\left(\overline{p}_{i} +\left(\frac{p_{i} +\left(2^{\alpha } -1\right)\overline{p}_{i}}{2^{\alpha } } \right)\right)-\left|\overline{p}_{i} -\left(\frac{p_{i} +\left(2^{\alpha } -1\right)\overline{p}_{i} }{2^{\alpha } } \right)\right|} \end{array}\right) =2\] 
\[\Rightarrow 4-2^{1-\alpha } \sum _{i=1}^{n}\left|p_{i} -\overline{p}_{i} \right| =4                                       \Rightarrow \sum _{i=1}^{n}\left|p_{i} -\overline{p}_{i} \right| =0\] 
\textit{i. e.} $p_i=\overline{p}_{i}$ for each \textit{i}. Conversely if $p_i=\overline{p}_{i}$ for each \textit{i}, then clearly${\ I}_{\alpha }\left(P(n);\overline{P}(n)\right)=0$. This completes the proof of \textbf{M1}. Further the proof of \textbf{M2 }is obvious since the measure ${\ I}_{\alpha }\left(P(n);\overline{P}(n)\right)$ is symmetric with respect to \textit{P} and \textit{Q.} For proof of \textbf{M3}, consider the following inequalities\textbf{}
\[min\left(p_{i} ,\overline{p}_{i} \right)\le min\left(p_{i} ,\left(\frac{p_{i} +\overline{p}_{i} }{2} \right)\right)\le min\left(p_{i} ,\left(\frac{3p_{i} +\overline{p}_{i} }{4} \right)\right)\le ...\] 
\[\le min\left(p_{i} ,\left(\frac{\left(2^{\alpha } -1\right)p_{i} +\overline{p}_{i} }{2^{\alpha } } \right)\right)\le ...\forall i\]
and
\[min\left(p_{i} ,\overline{p}_{i}\right)\le min\left(\left(\frac{p_{i} +\overline{p}_{i} }{2} \right),\overline{p}_{i} \right)\le min\left(\left(\frac{p_{i} +3\overline{p}_{i} }{4} \right),\overline{p}_{i} \right)\le...\] 
\[\le min\left(\left(\frac{p_{i} +\left(2^{\alpha } -1\right)\overline{p}_{i} }{2^{\alpha } } \right),\overline{p}_{i} \right)\le ...\forall i\]
This gives
\[min\left(p_{i} ,\overline{p}_{i} \right)\le \frac{1}{2} \left(min\left(p_{i} ,\left(\frac{p_{i} +\overline{p}_{i} }{2} \right)\right)+min\left(\left(\frac{p_{i} +\overline{p}_{i} }{2} \right),\overline{p}_{i} \right)\right)\]
\[\le \frac{1}{2} \left(min\left(p_{i} ,\left(\frac{3p_{i} +\overline{p}_{i} }{4} \right)\right)+ min\left(\left(\frac{p_{i} +3\overline{p}_{i} }{4} \right),\overline{p}_{i} \right)\right)\]
\[\le ...\le \frac{1}{2}\left(min\left(p_{i} ,\left(\frac{\left(2^{\alpha } -1\right)p_{i} +\overline{p}_{i} }{2^{\alpha } } \right)\right)+min\left(\left(\frac{p_{i} +\left(2^{\alpha } -1\right)\overline{p}_{i} }{2^{\alpha } } \right),\overline{p}_{i} \right)\right)\le ...\forall i\]
\[{\Rightarrow -log\left(\frac{1+\sum _{i=1}^{n}min\left(p_{i} ,\overline{p}_{i} \right) }{2} \right)\le -log\left(\frac{1+\frac{1}{2} \sum _{i=1}^{n}\left(min\left(p_{i} ,\left(\frac{p_{i} +\overline{p}_{i} }{2} \right)\right)+min\left(\left(\frac{p_{i} +\overline{p}_{i} }{2} \right),\overline{p}_{i} \right)\right) }{2} \right)} \]
\[{\le -log\left(\frac{1+\frac{1}{2} \sum _{i=1}^{n}\left(min\left(p_{i} ,\left(\frac{3p_{i} +\overline{p}_{i} }{4} \right)\right)+min\left(\left(\frac{p_{i} +3\overline{p}_{i} }{4} \right),\overline{p}_{i} \right)\right) }{2} \right)}\]
\[{\le ...\le -log\left(\frac{1+\frac{1}{2} \sum _{i=1}^{n}\left(min\left(p_{i} ,\left(\frac{\left(2^{\alpha } -1\right)p_{i} +\overline{p}_{i} }{2^{\alpha } } \right)\right)+min\left(\left(\frac{p_{i} +\left(2^{\alpha } -1\right)\overline{p}_{i} }{2^{\alpha } } \right),\overline{p}_{i} \right)\right) }{2} \right)\le ...\forall i} \] 
\[\Rightarrow I_{0} (P(n);\overline{P}(n))\le I_{1} (P(n);\overline{P}(n))\le I_{2} (P(n);\overline{P}(n))\le ...\le I_{\alpha } (P(n);\overline{P}(n))\le ...\forall i\]
This completes the proof of theorem 3.1. Also from (1.7) and (1.9), we can write 
\[ I_{\alpha}(P(n);\overline{P}(n))\le I_{\alpha} (P(n);\overline{\overline{P}}(n))\le I_{\alpha} (P(n);\overline{\overline{\overline{P}}}(n))\le ...\le ...\forall i\]
for all $\alpha =0,1,2,3,...$.
\section{Conclusion}
In the present work, we have established some new inequalities involving probability functions $P(n)$ and $\overline{P}(n)$. Also we have developed a parametric dissimilarity function which quantifies the dissimilarity between $P(n)$ and $\overline{P}(n)$. Work on the generalization of results obtained is in progress and will be communicated elsewhere.

\end{document}